\input amstex
\magnification =\magstep 1
\documentstyle{amsppt}
\pageheight{9truein}
\pagewidth{7truein}
\NoRunningHeads
\nologo
\baselineskip=16pt

\topmatter
\title On conjugacy classes in  groups
\endtitle

\author Marcel Herzog*, Patrizia Longobardi** and Mercede Maj**
\endauthor

\affil *School of Mathematical Sciences \\
       Tel-Aviv University \\
       Ramat-Aviv, Tel-Aviv, Israel
{}\\
       **Dipartimento di Matematica \\
       Universit\`a di Salerno\\
       via Giovanni Paolo II, 132, 84084 Fisciano (Salerno), Italy
\endaffil

\thanks This work was supported by the National Group for Algebraic and
Geometric Structures, and their Applications (GNSAGA - INDAM), Italy.
The first author  is grateful to the Department of
Mathematics of the University of Salerno
for its hospitality and
support, while this investigation was carried out.
\endthanks
\endtopmatter

\document

\heading I. Introduction\\
\endheading

Let  $G$ be a group and let $x$ be an element of $G$. Then clearly $\langle x\rangle \leq C_G(x)$.
Thus $x$ satisfies either $\langle x\rangle = C_G(x)$ or $\langle x\rangle <  C_G(x)$.
Accordingly, we introduce the following definitions.
\definition {Definitions} Let $G$ be a group. Then:

(a) $G^{*}=G\setminus \{1\}$.

(b) An element $x$ of $G^{*}$ will be called {\it deficient}  if
$$ \langle x\rangle < C_G(x)$$
and it will be called {\it non-deficient} if
$$\langle x\rangle = C_G(x).$$

If $x\in G$ is deficient (non-deficient), then the conjugacy class $x^G$ of $x$ in $G$ will be also called {\it deficient}
({\it non-deficient}).

(c) Let $j$ be a {\it non-negative} integer.
We shall say that the  group $G$ has {\it defect $j$}, denoted by $G\in D(j)$
or by the phrase "$G$ is a $D(j)$-group",
if {\bf exactly}  $j$  non-trivial conjugacy classes of $G$ are deficient.

(d) Moreover, if $G$ is a periodic group, then  $\pi(G)$ will denote the set  of all primes dividing the order of some element of $G$.
\enddefinition

In this paper a group is called {\it arbitrary} if it is either finite or infinite.

Our paper consists of three parts: Part A, dealing with finite $D(0)$-groups and $D(1)$-groups and consisting of Sections II-V, Part B,
dealing with some preliminary results from our paper [6] and consisting of Sections VI-VII,  and Part C,
dealing with arbitrary $D(0)$-groups and $D(1)$-groups and consisting of Sections VIII-IX.

The aim of Part A  is to determine all finite $D(0)$-groups and $D(1)$-groups.
We prove the following two theorems. In their statement, and elsewhere, we shall use the following notation. The phrase
"$G=A\rtimes B$ is a Frobenius group" means that $G$ is a finite Frobenius group with the kernel $A$
and a complement $B$. It is well known that $|A|\equiv 1\pmod{|B|}$.
Moreover, $C_n$ denotes the cyclic group of order $n$,
$D_n$ denotes the dihedral group of order $n$, $Q_8$ denotes the quaternion group
and if $p$ is a prime and $s$ is a positive integer, then $E_{p^s}$ denotes the elementary
abelian group of order $p^s$.

\proclaim {Theorem 1} Let $G$ be a finite group. Then $G\in D(0)$ if and only if one of the following statements holds:
\roster
\item $G$ is an abelian group and either $G=\{1\}$, or $G=C_p$  for some prime $p$,

\item $G$ is a non-abelian group  and $G=C_p\rtimes C_q$ is a Frobenius group,
with $p$ and $q$ distinct primes.
\endroster
\endproclaim
\proclaim {Theorem 2} Let $G$ be a finite group. Then $G\in D(1)$ if and only if one of the following statements  holds:
\roster
\item $G$ is a nilpotent group and either $G=C_4$ or $G=Q_8$,
\item $G$ is a non-nilpotent solvable group,  and either  $G=E_{2^s}\rtimes C_q$
is a Frobenius group with  a positive integer $s$ and $q=2^s-1$ being a Mersenne prime,
or $G=C_q\rtimes C_4$ is a Frobenius group with $q$ being a prime satisfying $q\equiv 1\pmod{4}$,
or $G=D_{18}$,
\item $G$ is a non-solvable group and either $G=A_5$ or $G=PSL(2,7)$.
\endroster
\endproclaim

In Part B we shall quote some results from our paper [6], dealing with arbitrary groups $G$, such that each
element of $G^{*}$ is of prime power order and $G$ satisfies
certain additional conditions. These results will be very useful in our Part C.

Part C deals with arbitrary $D(0)$-groups and $D(1)$-groups.
Our aim in Part C is to find properties of arbitrary $D(0)$-groups and $D(1)$-groups, which force these groups to be finite.

Concerning the $D(0)$-groups we prove the following  four theorems.

\proclaim {Theorem 3} Let $G$ be a  locally finite $D(0)$-group. Then $G$ is a finite group.
\endproclaim

\proclaim {Theorem 4} Let $G$ be a  locally graded $D(0)$-group. Then $G$ is a finite group.
\endproclaim

\proclaim {Theorem 5} Let $G$ be a residually finite $D(0)$-group. Then $G$ is a  finite group.
\endproclaim

\proclaim {Theorem 6}  Let $G$ be a $D(0)$-group and suppose that $2\in \pi(G)$. Then $G$ is a finite group.
\endproclaim

Finally, concerning  the $D(1)$-groups we prove the following  three theorems.

\proclaim{Theorem 7} Let $G$ be a locally finite $D(1)$-group. Then $G$ is a finite group.
\endproclaim

\proclaim{Theorem 8} Let $G$ be a periodic locally graded $D(1)$-group. Then $G$ is a finite group.
\endproclaim

\proclaim{Theorem 9} Let $G$ is a residually finite $D(1)$-group. Then $G$ is a finite group.
\endproclaim

Theorem 1 is proved in Section II, Theorem 2 is proved in Section V, Theorems 3, 4, 5 and 6 are proved in Section VIII
and Theorems 7, 8 and 9 are proved in Section IX.

Our Introduction is now complete.

\heading Part A - Finite $D(0)$-groups and $D(1)$-groups\\
\endheading

The aim of Part A is to prove Theorem 1 and Theorem 2 stated in the Introduction.
Theorem 1 will be proved in Section II, and Theorem 2 will be proved in Section V.

\heading II. Finite $D(0)$-groups and a proof of Theorem 1\\
\endheading
The aim of this section is to prove Theorem 1, which determines all finite $D(0)$-groups.
For our proof we need the following two preliminary results.
The first one is the following proposition.

\proclaim{Proposition A1} The finite group $G$ satisfies $G\in D(0)$ if and only if the following
two statements hold:
\roster
\item All elements of $G^{*}$ are of prime order, and
\item If $p\in \pi(G)$, then a Sylow $p$-subgroup of $G$ is of prime order.
\endroster
\endproclaim
\demo{Proof} Suppose, first, that (1) and (2) hold. It follows by (1) that each $x\in G^{*}$ is
of prime order, say $p$, and (2) implies that $\langle x\rangle$ is a Sylow $p$-subgroup of $G$. If
$y\in C_G(x)^{*}$, then $|y|=q$ for some prime $q$. If $q\neq p$,   then $xy\neq 1$ and it is not of prime order,
in contradiction to (1). Hence if $y\in C_G(x)^{*}$, then $|y|=p$.
If $y\notin \langle x\rangle$, then $\langle x,y\rangle$ is a $p$-group
of order $p^2$, in contradiction to (2). Hence $C_G(x)=\langle x\rangle$ for all $x\in G^{*}$
and $G\in D(0)$, as required.

Conversely, if $C_G(x)=\langle x \rangle $ for all $x\in G^{*}$, then clearly all elements of
$G^{*}$ are of prime order and for each  $p\in \pi(G)$,  a Sylow $p$-subgroup
of $G$ is of prime order. Hence $G$ satisfies (1) and (2), as required.
\qed
\enddemo

We shall denote the statement (i) of Proposition A1 by (A1(i)).
A similar notation will be used concerning our other results.

The second preliminary result is the following theorem of
K.N. Cheng, M. Deaconescu, M.L. Lang and W. Shi (see [1]).
\proclaim {Theorem N} Suppose that $G$ is a finite group and all elements of $G^{*}$ are of prime order. Then one of the
following statements holds.
\roster
\item $G$ is a nilpotent group if and only if either $G=\{1\}$ or $G$ is a $p$-group of exponent $p$ for
some prime $p$.
\item $G$ is a non-nilpotent solvable group if and only if $G=P\rtimes C_q$ is a Frobenius
group, where $P$ is the Sylow $p$-subgroup of $G$ of exponent $p$ for some prime $p$
and $q$ is a prime different from $p$.
\item $G$ is a non-solvable group if and only if $G=A_5$.
\endroster
\endproclaim

We are now ready for the proof of  Theorem 1.
\demo{Proof of Theorem 1} By Proposition A1, all groups mentioned in Theorem 1 belong to $D(0)$.

Conversely, if $G\neq \{1\}$ is a finite $D(0)$-group, then it satisfies statements (1) and (2) of
Proposition A1. Groups satisfying (A1(1)) are described in Theorem N. Since $G$ also satisfies (A1(2)),
statement (3) of Theorem N is impossible and in statement (2) of Theorem N we must have
$|P|=p$ and $G=C_p\rtimes C_q$ is a Frobenius group
with $p$ and $q$ distinct primes.
Finally, by  (A1(2)), in statement (1) of Theorem N
we must have $G=C_p$. Thus either (1) or (2) of Theorem 1 holds, as required.
\qed
\enddemo

\heading III. Properties of finite $D(1)$-groups \\
\endheading

We now turn to finite $D(1)$-groups. We begin with two propositions dealing with properties of such groups.
\proclaim{Proposition A2}  Let $G$ be a finite group. If $G\in D(1)$, then the following statements hold.
\roster
\item All elements of $G^{*}$ are either of prime order or of prime-squared order.
\item There exists $p\in \pi(G)$ such that a Sylow $p$-subgroup $P$
of $G$ is of
order $|P|=p^s\geq p^2$, and for all other $q\in \pi(G)$, a Sylow $q$-subgroup
of $G$ is of the prime order $q$. In particular, either $exp(P)=p$ or $exp(P)=p^2$.
\item All elements of $G$ of order $p$ are conjugate in $G$.
\item If $exp(P)=p^2$, then either $P=C_{p^2}$ or $|Z(P)|=p$ and $|P|=p^s\geq p^3$.
\item If $|P|=p^s\geq p^3$, then either $exp(P)=p$ or $exp(P)=p^2$ and $|Z(P)|=p$.
In any case, $Z(P)$ is an elementary abelian $p$-group.
\item If $p>2$, then $G$ is solvable.
\item If $G> P$, then $Z(G)=\{1\}$.
\item If $R$ is a non-trivial normal Sylow $r$-subgroup of $G$ for some $r\in \pi(G)$, then $G$ is a Frobenius group with the kernel $R$.
\endroster
\endproclaim
\demo{Proof} \roster
\item Suppose that $x\in G$ and $rt\bigm| |x|$, where $r$ and $t$ are distinct primes.
Let $y=x^r$ and $z=x^t$. Then $|y|\neq |z|$, $|y|<|x|$, $|z|<|x|$ and $x\in C_G(y)\cap C_G(z)$.
Hence
$y$ and $z$ are
elements of $G$ of different orders and both are deficient,
in contradiction to $G\in D(1)$. Therefore all elements  of $G^{*}$ are of prime-power
order. If $x\in G$ and $|x|=r^3$ for some prime $r$, then $x^r$ and $x^{r^2}$ are
elements of $G$ of different orders and both are deficient, contradicting $G\in D(1)$. Therefore
statement (1) holds.
\item Suppose that $P\in Syl_p(G)$ and $Q\in Syl_q(G)$ for distinct prime divisors $p$ and $q$ of $|G|$
and both $P$ and $Q$ are not of prime order. Then their centers contain a deficient element
of prime order, contradicting $G\in D(1)$. So a Sylow subgroup of $G$ is not of prime order for
at most one prime divisor of $|G|$. If all subgroups of $G$ are of prime order, then in view of
statement (1), Proposition A1 implies that $G\in D(0)$, contradicting $G\in D(1)$. Therefore
there exists $p\in \pi(G)$ such that a Sylow $p$-subgroup $P$
of $G$ is of
order $|P|=p^s\geq p^2$, and for all other $q\in \pi(G)$, a Sylow $q$-subgroup
is of the prime order $q$, as required. Moreover, by (1) either $exp(P)=p$ or $exp(P)=p^2$.
\item Since $|P|\geq p^2$, all elements of $G$ of order $p$ are deficient, and since $G\in D(1)$,
all elements of $G$ of order $p$ are conjugate in $G$, as required.
\item If $exp(P)=p^2$, then there exists $x\in P$ of order $p^2$. Since $x^p$ is a deficient element
of $G$, $G\in D(1)$ implies $x$ is a non-deficient element of $G$ and hence $Z(P)\leq C_P(x)=\langle x\rangle$.
If $x\in Z(P)$, then $P=C_P(x)=\langle x\rangle$ and $P=C_{p^2}$. If $x\notin Z(P)$,
then $1<Z(P)< C_P(x)=\langle x\rangle$, which implies that $Z(P)=\langle x^p\rangle$ is of order $p$.
Hence $P$ is non-abelian, which implies that $|P|=p^s\geq p^3$, as required.
\item Suppose that $|P|\geq p^3$. Then, by (2), either $exp(P)=p$ or  $exp(P)=p^2$. If $exp(P)=p^2$, then $|Z(P)|=p$, by (4).
So in any case $Z(P)$ is an elementary abelian $p$-group.
\item If $p>2$, then the order of a Sylow $2$-subgroup of $G$ is either $1$ or $2$ and $G$ is solvable
by the Feit-Thompson theorem.
\item If $G> P$, then $|G|$ is divisible by at least two primes. If $ Z(G)\neq \{1\}$, then there exists $y\in Z(G)$ of prime order and hence
$G$ contains an element of order divisible by two primes, in contradiction to statement (1).
Hence if $G>P$, then $Z(G)=\{1\}$.
\item Since $R$ contains all $r$-elements of $G$, it follows by (1) that $C_G(g)\leq R$ for all $g\in R^{*}$ and the claim follows.
\qed
\endroster
\enddemo

\remark {Remark} If $G$ is a finite $D(1)$-group, then $p$ and $P$ will denote the prime $p$ and the Sylow $p$-subgroup  of $G$,  as defined in (A2(2)).
\endremark

We conclude this section with our first classification result concerning finite $D(1)$-groups.

\proclaim {Proposition A3} If  $G$ is a finite $D(1)$-group and $G=P$, then $p=2$
and either $G=C_4$ or $G=Q_8$, and vice-versa.
\endproclaim
\demo{Proof} If $G=P$, then $|G|=p^s\geq p^2$. If either
$p>2$ or $p=2$ and $G=E_{2^2}$, then $Z(G)$ contains two distinct
deficient elements of order $p$, which are not conjugate to each other,
contradicting $G\in D(1)$. Hence $p=2$,  $|G|\geq 4$ and if $|G|=4$, then
$G=C_4$, which satisfies $G\in D(1)$.

Suppose, now, that $|G|=2^s\geq 8$. By (A2(3)), all involutions
are conjugate in $G$ and hence they  belong to $Z(G)$. It follows that
$G$ contains only one involution and either $G=C_{2^s}$ or $G$ is a generalized quaternion group of order
$2^s$. Since by (A2(2)) $\exp(G)\leq 4$, $G=C_{2^s}$ is impossible and so is $Q_{2^s}$ for $2^s\geq 16$.
Hence $G=Q_8$, which satisfies $G\in D(1)$.

Thus either $G=C_4$ or $G=Q_8$, and vice-versa, as required.
\qed
\enddemo

\heading IV. Classification of finite solvable $D(1)$-groups
\endheading

Finite solvable  groups satisfying  condition (A2(1)) were studied by Graham Higman
in [7]. He proved the following more general result.
\proclaim {Theorem GH} Let $G$ be a finite solvable group with all elements of prime-power order.
Let $r$ be a prime such that $G$ has a normal $r$-subgroup of order greater than $1$, and let $R$
be the greatest normal $r$-subgroup of $G$. Then $G/R$ satisfies one of the following
statements.
\roster
\item $G/R$ is a cyclic group, whose order is a power of a prime other than $r$.
\item $r$ is odd and $G/R$ is a generalized quaternion group.
\item $G/R$ is a group of order $r^at^b$, where  $t$ denotes a
prime of the form $kr^a+1$ and the Sylow subgroups of $G/R$ are cyclic.
\endroster
Thus $G$ has order divisible by at most two primes, and $G/R$ is metabelian.
\endproclaim

Let now $G$ be a finite solvable $D(1)$-group. In such case we shall say that $G\in SD(1)$
or "$G$ is an $SD(1)$-group".
Recall that the order of $G$ is divisible by a prime $p$, such that $P\in Syl_{p}(G)$
is not of  prime order, while all Sylow subgroups corresponding to other primes dividing $|G|$
are of prime order. In particular, $|P|=p^s\geq p^2$.

Since $G$ satisfies (A2(1)), it follows by Theorem GH that $|G|=p^sq^v$, where $q$ is a prime
different from $p$, $s\geq 2$ and $v\in \{0,1\}$. By Proposition A3, if $v=0$ then either
$G=C_4$ or $G=Q_8$.

So suppose, from now on, that $v=1$ and $|G|=p^sq$, and let $Q$ denote a Sylow $q$-subgroup
of $G$ of order $q$. Thus:
$$G=PQ, \quad\text{where}\ |P|=p^s\geq p^2\ \text{and}\ |Q|=q.$$
By (A2(2)), either $\exp(P)=p$  or $\exp(P)=p^2$. Recall that the normal subgroup $R$ of $G$
in Theorem GH satisfies one of the following conditions: either
$R=Q\triangleleft G$, or $R=P\triangleleft G$, or $R$ satisfies $1<R<P$ and it is
the maximal normal $p$-subgroup of $G$. Since $G\in D(1)$, these three cases are
disjoint.

In Proposition A3 all finite $D(1)$-groups $G$ satisfying $G=P$ were determined. In this section we shall complete
the classification of  finite solvable $D(1)$-groups. In Proposition A5, we shall classify the  $SD(1)$-groups $G$ satisfying
$G>P$ and $exp(P)=p$ and in  Proposition A6 we shall classify the   $SD(1)$-groups $G$ satisfying
$G>P$ and $exp(P)=p^2$. Finally, in Proposition A7, we shall present the classification of all finite solvable $D(1)$-groups.

But first, we need to prove the following proposition. We shall use the following notation:
$$\Omega_1(P)=\langle x| x\in P\ \text{and}\ |x|=p\rangle.$$

\proclaim {Proposition A4} Let $G\in SD(1)$. Then
 the following statements hold.
 \roster
 \item If $P\neq P^u$ for some $u\in G$, then
 $$Z(P)\cap Z(P^u)=1.$$
 \item If $p>2$, then $P\triangleleft G$.
 \item If $P\ntriangleleft G$, then $p=2$.
 \item If $P\triangleleft G$, then $\Omega_1(P)\leq Z(P)$ and $\exp(\Omega_1(P))=p$.
 \item If $Q\triangleleft G$, then $P$ is cyclic.
 \endroster
 \endproclaim

 \demo{Proof} \roster
 \item If $P\neq P^u$ for some $u\in G$, then $G>P$. If $x\in Z(P)\cap Z(P^u)$, then $C_G(x)>P$
 and since $G=PQ$, it follows that $x\in Z(G)$ and by (A2(7)) $x=1$, as required.
 \item Since $p>2$, there exist $u,v\in Z(P)$ with $u\neq v$ and $|u|=|v|=p$. By (A2(3)),
 $u$ and $v$ are conjugate in $G$ and there exists $z\in Q^{*}$ such that $v=u^z$. If
 $z\in N_G(P)$, then $Q \leq N_G(P)$ and $P\triangleleft G$, as required. Suppose, finally, that
 $z\notin N_G(P)$. Then $P^z\neq P$ and by (1) $Z(P)\cap Z(P^z)=1$. But $u,v\in Z(P)$ and
 $v=u^z\in Z(P^z)$, so $v\in Z(P)\cap Z(P^z)$, a contradiction. Therefore $P\triangleleft G$, as
 required.
 \item Follows  by (2).
 \item By (A2(3)),  all elements of order $p$ in $P$ are conjugate to each other in $G$. Since
 $P\triangleleft G$, also $Z(P)\triangleleft G$ and hence all elements of $P$ of order $p$ belong to $Z(P)$. Thus $\Omega_1(P)\leq Z(P)$
 and $\exp(\Omega_1(P))=p$, as claimed.
 \item If $Q\triangleleft G$, then it follows by (A2(8)) that $G=Q\rtimes P$ is a Frobenius group.
 Since $C_G(Q)=Q$ and $N_G(Q)=G$,
 it follows that $P$ is isomorphic to $N_G(Q)/C_G(Q)$ and since $Q$ is cyclic of prime order,
 also $P$ is cyclic, as claimed.
 \qed
 \endroster
 \enddemo

 We proceed with the proof of our first main result of this section.
 \proclaim {Proposition A5} If $G\in SD(1)$, $G>P$ and $\exp(P)=p$,
 then $p=2$, $P=E_{2^s}$, $P\triangleleft G$, $q=2^s-1 $ is a Mersenne prime and
 $$G=E_{2^s}\rtimes C_q$$
is a Frobenius group, and vice-versa.
\endproclaim

\demo{Proof} Suppose that $G\in SD(1)$, $G>P$ and $\exp(P)=p$.

Let $z\in Z(P)^{*}$. Since $\exp(P)=p$, it follows that $|z|=p$ and by (A2(3)) all elements of order $p$ in $G$ are conjugate.
Hence $z^G=z^Q$ contains all elements of $G$ of order $p$ and each element of $z^G=z^Q$ belongs to the center of one of the Sylow $p$-subgroups of $G$.
Since  the non-unit
elements of each such Sylow subgroup belong to $z^G$, it follows that
$z^G=z^Q=\cup_{P^u} Z(P^u)^{*}$, where the union
extends over all the distinct Sylow $p$-subgroups $P^u$ of $G$.

If $P\ntriangleleft G$, then $[G:N_G(P)]=q$ and  $|z^Q|=q$.
Since by (A4(1)) $Z(P^u)\cap Z(P^v)=1 $ for all $P^u\neq P^v$, it follows that $q=|z^Q| =q(|Z(P)|-1)$, implying that $|Z(P)|=2$.
Hence $p=2$, and since $exp(P)=p$, it follows that $P$ is abelian. But then $2=|Z(P)|=|P|$, a contradiction, since $|P|=2^s\geq 2^2$.

So $P\triangleleft G$, $P^{*}= z^Q$ and $q=|z^Q|=p^s-1=(p-1)(p^{s-1}+\dots +1)$, which implies that $p=2$ and
$exp(P)=2$. Hence $P$ is an elementary abelian $2$-group and $q=2^s-1$ is a Mersenne prime. Thus $G=E_{2^s}\rtimes C_q$
and by (A2(8)) $G$ is a Frobenius group, as required.

Conversely, we need to show that if $q=2^s-1$ is a Mersenne prime, then there exists a Frobenius group  $G=E_{2^s}\rtimes C_q$, which
belongs to $D(1)$.

First notice that $Aut(E_{2^s})=GL(s,2)$ and as shown in Satz 7.3 on page 187 of Huppert's book [8],
$GL(s,2)$  contains a cyclic subgroup $T$ of order $2^s-1$ (called the Singer-cycle), which
acts regularly on $E_{2^s}$. Hence if $q=2^s-1$ is a Mersenne prime, then  a Frobenius group $G=E_{2^s}\rtimes C_q$ exists.
Since all involutions in $G$ are conjugate and
all other non-trivial elements are of order $q$, with centralizers of order $q$, it follows that $G$ is
a $D(1)$-group, as required.

The proof of Proposition A5 is now complete.
\qed
\enddemo

We proceed with the proof of our second main result of this section.
\proclaim {Proposition A6} If $G\in SD(1)$, $G>P$ and $\exp(P)=p^2$, then
$P=C_{p^2}$ and one of the following statements holds:
\roster
\item $P\triangleleft G$, $p=3$, $q=2$ and $G=D_{18}$,
\item $Q\triangleleft G$, $p=2$, $q\equiv 1\pmod{4}$ and $G=C_q\rtimes C_4$, a Frobenius group,
\endroster
and vice-versa.
\endproclaim

\demo{Proof} Suppose that $G\in SD(1)$, $G>P$ and $\exp(P)=p^2$. It follows by (A2(4)) that
either $P=C_{p^2}$ or $|Z(P)|=p$ and  $|P|=p^s\geq p^3$.

Suppose, first, that $|Z(P)|=p$ and  $|P|=p^s\geq p^3$. Then $P$ is non-abelian.
Moreover, all elements of $G$ of order $p$ are clearly deficient, so they are conjugate to each other.

If $P\triangleleft G$, then $Z(P)\triangleleft G$ and since $|Z(P)|=p$, $Z(P)^{*}=z^Q$
for each $z\in Z(P)^{*}$. Hence $p-1=q$, implying that $q=2$. Since $P\triangleleft G$, it follows by (A2(8))
that $G=P\rtimes C_2$ is a Frobenius group and hence $P$ is abelian, a contradiction.

Therefore $P\ntriangleleft G$ and $p=2$ by (A4(3)). Hence $|Z(P)|=2$ and by (A2(1))
the elements of $P^{*}$ are either involutions or elements of order $4$. Moreover, since $\exp(P)=p^2$,  $P$ contains an element of
order $4$. If $x\in P$ is of order
$4$, then $x^2$ is deficient and hence $x$ is non-deficient. It follows that the deficient elements of $G$ are
of order $2$ and $Z(P)<\langle x\rangle$, implying that $Z(P)=\langle x^2\rangle$. Since that is true for each element of $P$ of
order $4$, it follows that
$P/Z(P)$ is of exponent $2$. Hence $P/Z(P)$ is an abelian group and $P'=Z(P)$. Since $|P'|=|Z(P)|=2$,
it follows that $P$ is an extra-special $2$-group. Thus $P$ is a central-product of $T_1,T_2,\dots, T_n$,
with each $T_i$ being equal  either to $D_8$ or to $Q_8$. If $n\geq 2$, then $T_2\leq C_P(T_1)$
and the elements of order $4$ in $T_1$ are deficient, contrary to $G\in SD(1)$. Hence $P$
is either $D_8$ or $Q_8$. Since $P\ntriangleleft G$ and by (A4(5)) also $Q\ntriangleleft G$, it follows
by the GH-theorem that there exists $R<G$ satisfying $1<R<P$ and $R\triangleleft G$. If $|R|=2$ or if
$R=C_4$, then $Z(G)>1$, in contradiction to (A2(7)). Hence we must have $R=E_{2^2}$ and $P=D_8$.
Let $v$ be an involution in $P$ outside $R$. Then $v$ is not conjugate to an involution in $R$,
a contradiction.

We have shown that $|P|=p^s\geq p^3$ is impossible. Hence
$P=C_{p^2}$ holds.

If $Q\triangleleft G$, then $P\ntriangleleft G$ and by (A4(3)) $p=2$. Hence $G=Q\rtimes C_4$.
By (A2(8)) $G$ is a  Frobenius group, implying that $q\equiv 1\pmod{4}$, and (2) holds.

Conversely, we need to show that if $q$ is a prime number satisfying $q\equiv 1\pmod{4}$, then
there exists a Frobenius group $G=Q\rtimes C_{4}$, with $Q$ being a group of order $q$ and $G$ belonging to $D(1)$.
Indeed, if $Q$ is a group of order $q$, then $q\equiv 1\pmod{4}$ implies that $C_4$
is a subgroup of $Aut(Q)$, so there exists a group $G=Q\rtimes C_{4}$, with $C_4\leq Aut(Q)$. Since $Q$ is
of prime order, this group is a Frobenius group with the kernel $Q$ and with the complements conjugate to $C_4$.
Thus $G$ has one conjugacy class of involutions. Since
only the involutions in $G$ are deficient elements, it follows that $G$ belongs to $D(1)$, as required.

Suppose, now, that $P=C_{p^2}\triangleleft G$. Then $\Omega_1(P)\triangleleft G$, $|\Omega_1(P)|=p$ and
$\Omega_1(P)^{*}=z^Q$ for each $z\in  \Omega_1(P)^{*}$. Hence $p-1=q$, which implies that $p=3$
and $q=2$. It follows by (A2(8)) that $G=C_9\rtimes C_2$ is a Frobenius group, so $G=D_{18}$ and (1) holds.
Conversely, it is easy to see that $D_{18}\in SD(1)$, as required.

Suppose, finally, that $Q\ntriangleleft G$ and $P \ntriangleleft G$. Since $P=C_{p^2}$, the GH-theorem
implies that $C_p\triangleleft G$. Then $C_p\rtimes Q$ is a non-abelian group of order $pq$, with
$q\mid p-1$. Thus $q<p$ and therefore $P\triangleleft G$, a contradiction.

The proof of Proposition A6 is now complete.
\qed
\enddemo

Propositions A3, A5 and A6 imply the following final result.
\proclaim {Proposition A7} Let $G$ be a finite solvable group and let $P$ be a Sylow $p$-subgroup of $G$ for a prime $p$. Then $G\in D(1)$ if and only if one of the following cases holds.
\roster
\item $G=P$, $p=2$ and either $G=C_4$ or $G=Q_8$.
\item $G>P$, $p=2$, $P=E_{2^s}$ with $q=2^s-1$ being a Mersenne prime and
$$G=E_{2^s}\rtimes C_q$$
is a Frobenius group.
\item $G>P$, $p=3$, $P=C_{3^2}$, $q=2$ and $G=D_{18}$.
\item $G>P$, $p=2$, $q$ is a prime satisfying $q\equiv 1\pmod {4}$, $P=C_4$ and
$$G=C_q\rtimes C_4$$
is a Frobenius group.
\endroster
\endproclaim
\demo {Proof} Case (1) follows from Proposition A3, Case (2) follows from Proposition A5 and Cases (3) and (4) follow from Proposition A6.
\qed
\enddemo

\heading V. Finite non-solvable  $D(1)$-groups and a proof of Theorem 2\\
\endheading

Our first aim in this section is to determine all finite non-solvable  $D(1)$-groups.
We start with a series of definitions.
\definition{Definitions}
\roster
\item A group $G$ is a $CIT$-group if $G$ is a finite group of even order and the centralizer of
every involution in $G$ is a $2$-group.
\item A simple group means a non-abelian simple group.
\item A group $G$ is an $NSD(1)$-group, or $G\in NSD(1)$, if $G$ is a finite non-solvable  $D(1)$-group.
\endroster
\enddefinition

The non-solvable finite $D(1)$-groups are determined in the following proposition.
\proclaim{Proposition A8} If $G$ is a $ NSD(1)$-group, then either $G=A_5$ or $G=PSL(2,7)$,
and vice-versa.
\endproclaim
\demo{Proof} It is easy to check that the simple groups $A_5$ and $PSL(2,7)$
are $D(1)$-groups. Therefore it suffices to prove that if
$G\in NSD(1)$, then either $G=A_5$ or $G=PSL(2,7)$.

So suppose that $G\in NSD(1)$. Then by the Feit-Thompson theorem the order of $G$ is divisible by $4$.
Hence by (A2(2)) $p=2$ and by (A2(3)) all involutions of $G$ are conjugate in $G$.
Moreover, by (A2(1)) the centralizer of
every involution in $G$ is a $2$-group. Thus $G$ is a non-solvable
$CIT$-group. By Theorem 5 in Suzuki's paper [10], the maximal solvable normal subgroup $N$ of $G$
is a $2$-group, and $G/N$ is either a simple $CIT$-group or $M_9$ in the notation of
Zassenhaus. As indicated by Suzuki, the group $M_9$ contains $PSL(2,9)$ as a
normal subgroup of index $2$. Thus $G/N$ is of even order.
Moreover, by Theorem 5 on page 438 of [11], if $G/N$ is of even order,
then $G$ contains
an involution outside $N$. Since all involutions in $G$
are conjugate to each other, it follows that $N=\{1\}$ and $G$  is either
a simple $CIT$-group or $M_9$. Since by (A2(1)), if $G$ is a finite $D(1)$-group, then all elements of $G^{*}$ are of prime-power
order, it follows  by Theorem  16 in [12] that simple finite $D(1)$-groups are either $PSL(2,q)$ for
$q=5,7,8,9,17$, or one of the groups $PSL(3,4)$, $Sz(8)$ and $Sz(32)$. Hence $G$ is either one of the above
simple groups or $M_9$. Since $M_9$ contains the group $PSL(2,9)$ as a normal subgroup of index $2$
and $PSL(2,9)$ has deficient elements of orders $2$ and $3$,  $M_9$ is not a $D(1)$-group.
It is also easy to check  that
$PSL(2,q)$ for $q=8,9,17$, $PSL(3,4)$, $Sz(8)$ and $Sz(32)$ are not $D(1)$-groups. Hence the only
remaining candidates are $A_5$ and $PSL(2,7)$, as claimed. The proof of Proposition A8 is now complete.
\qed
\enddemo

Our final aim in this section is to prove Theorem 2.
\demo{Proof of Theorem 2} By Propositions A7 and A8, all groups mentioned in Theorem 2 belong to $D(1)$.

Conversely,  suppose that $G$ is a finite $D(1)$-group.
If $G$ is solvable, then by Proposition A7 one of the statements (1) or (2) of Theorem 2 holds. If $G$ is non-solvable,
then by Proposition A8 the statement (3) of Theorem 2 holds. Hence the proof of Theorem 2 is now complete.
\qed
\enddemo

This is the end of our Part A.
Before moving to Part C, which  deals with $D(j)$-groups that are either finite or infinite,
we shall present some  results which were obtained in our paper [6]. These results are of independent interest and
they will be very helpful in Part C.

\heading  Part B - $CP$-groups satisfying some boundedness conditions\\
\endheading

In Part B, a group $G$ is either finite or infinite.
If $G$ is a periodic group, then  $\pi(G)$ will denote the set  of all primes dividing the order of some element of $G$.

Following   A.L. Delgado and Y. Wu in their paper [2], groups with each non-trivial element of prime power order will be called  $CP$-groups. Such groups are of course
periodic. In Part B we shall deal with $CP$-groups, which satisfy some boundedness condition, as defined below.

\definition {Definitions} A group $G$ will be  called a $BCP$-group if each element of $G^{*}$ is of prime power order and  for each $p\in \pi(G)$ there exists a positive integer
$u_p$  such that each $p$-element of $G$ is of order $p^i\leq p^{u_p}$.

A group $G$ will be called a $BSP$-group if each element of $G^{*}$ is of prime power order and for each $p\in \pi(G)$ there exists a positive integer
$v_p$  such that each finite $p$-subgroup of $G$ is of order $p^j\leq p^{v_p}$.
\enddefinition

The $BCP$-groups and the $BSP$-groups are clearly periodic. Notice that each $BSP$-group is a  $BCP$-group and each
$BCP$-group is a $CP$-group. Moreover, the $BCP$-property   and the $CP$-property are  inherited by subgroups and quotient
groups, and hence by sections. The $BSP$-
property is inherited by subgroups.

First we present our  results in [6] concerning $BCP$-groups.

\heading VI. Properties of $BCP$-groups\\
\endheading

It is well known  that  finitely generated groups  have only a
finite number of subgroups of a {\it given} finite index.
In particular, each such group has only a finite number of normal subgroups of a {\it given} finite index. We proved, using
the Zelmanov positive solution of the Restricted Burnside Problem (see [13] and [14])
that
finitely generated $BCP$-groups have only a finite number of normal subgroups of  {\it any} finite index.
\proclaim {Theorem B1} Let $G$ be a finitely generated $BCP$-group. Then $G$ has only a finite number of normal
subgroups of finite index.
\endproclaim

This basic result was used in the proofs of the next theorems dealing with $BCP$-groups.

Recall that a group $G$ is residually
finite if for each non-trivial element $g\in G$  there exists a normal subgroup $M(g)$ of $G$ such that $g\notin M(g)$ and $G/M(g)$ is finite.

As a corollary of Theorem B1  we obtained the following result.
\proclaim {Theorem B2} Let $G$ be a finitely generated residually finite $BCP$-group. Then $G$ is a finite group.
\endproclaim

It is well known that the residually finite property is inherited by subgroups. Hence Theorem B2 implies the following result.

\proclaim {Theorem B3} Let $G$ be a residually finite $BCP$-group. Then $G$ is a locally finite group.
\endproclaim

Recall that a group $G$ is locally graded if  each non-trivial finitely generated subgroup of $G$ has a proper normal subgroup of finite index.
By applying the above results, we proved the following theorem concerning  locally graded $BCP$-groups.

\proclaim {Theorem B4} Let $G$ be a locally graded $BCP$-group. Then $G$ is a locally finite group.
\endproclaim

As a corollary we get the following theorem.
\proclaim {Theorem B5} Let $G$ be a finitely generated locally graded $BCP$-group. Then $G$ is a finite group.
\endproclaim

Now we move to the properties of $BSP$-groups proved in [6].

\heading VII. Properties of $BSP$-gropups\\
\endheading
Recall that a group $G$ is a $BSP$-group if each element of $G^{*}$ is of prime power order and for each $p\in \pi(G)$ there exists a positive integer
$v_p$  such that each {\it finite $p$-subgroup} of $G$ is of order $p^j\leq p^{v_p}$.
This property is stronger than the $BCP$-property, which requires only that each {\it $p$-element} of $G$ is of order $p^i\leq p^{u_p}$.
Consequently, our results concerning the $BSP$-groups are stronger than those obtained for the $BCP$-groups.

The basic result concerning the $BSP$-groups is the following theorem.
\proclaim {Theorem B6} Let $G$ be a locally finite $BSP$-group. Then $G$ is a finite group.
\endproclaim
This theorem does not hold for $BCP$-groups, since if $p$ is a prime, then an infinite abelian $p$-group of finite exponent is
a locally finite $BCP$-group.

Theorem B6 yields the following strengthening of Theorem B4 for $BSP$-groups.

\proclaim {Theorem B7} Let $G$ be a  locally graded $BSP$-group. Then $G$ is a finite group.
\endproclaim

\demo {Proof} By Theorem B4 applied to $BSP$-groups, $G$ is a locally finite $BSP$-group. Hence, by Theorem B6, $G$ is a finite group, as required.
\qed
\enddemo

Finally, we also proved the following theorem.
\proclaim {Theorem B8} Let $G$ be a $BSP$-group and suppose that $2\in \pi(G)$. Then $G$ is a finite group.
\endproclaim

This is the end of Part B. We turn now to Part C, which is the last part of this paper.

\heading Part C - Arbitrary $D(0)$-groups and $D(1)$-groups\\
\endheading

A group which is either finite or infinite will be called {\it arbitrary}.
Part C deals with arbitrary $D(0)$-groups and $D(1)$-groups.

There exist infinite $D(j)$-groups.
In particular, the Tarski infinite $p$-groups, whose proper non-trivial subgroups are all of the prime order $p$,  are infinite $D(0)$-groups.
Our aim in Part C is to find properties of arbitrary $D(0)$-groups and $D(1)$-groups, which force these groups to be finite.

First we shall deal with arbitrary $D(0)$-groups.

\heading VIII. Arbitrary $D(0)$-groups and proofs of Theorems 3,4,5 and 6\\
\endheading

The above mentioned  Tarski infinite $p$-groups are infinite, periodic, non-locally-finite $D(0)$-groups.
Infinite  $D(0)$-groups were first studied by Delizia, Jezernik,
Moravec and Nicotera in their paper [3] (see also [4] and [5]).
In [3] they noticed that the following two theorems hold.

\proclaim {Theorem 3} Let $G$ be a  locally finite $D(0)$-group. Then $G$ is a finite group.
\endproclaim

\proclaim {Theorem 4} Let $G$ be a  locally graded $D(0)$-group. Then $G$ is a finite group.
\endproclaim

Using the results in Part B we shall provide different proofs of these theorems and also  proofs of the following two theorems.

\proclaim {Theorem 5} Let $G$ be a residually finite $D(0)$-group. Then $G$ is a  finite group.
\endproclaim

\proclaim {Theorem 6}  Let $G$ be a $D(0)$-group and suppose that $2\in \pi(G)$. Then $G$ is a finite group.
\endproclaim

\demo{Proofs of Theorems 3,4,5 and 6}
Our proofs of these theorems are based on the following observation.

\proclaim {Observation} An arbitrary $D(0)$-group is a $BSP$-group (and in particular a $BCP$-group).
\endproclaim

\demo{Proof} Indeed, if $G$ is an arbitrary $D(0)$-group, then $C_G(g)=\langle g\rangle$ for
each $g\in G^{*}$. If $x\in G$ is of infinite order, then
$$\langle x^2\rangle < \langle x\rangle \leq C_G(x^2),$$
a contradiction. Therefore each $D(0)$-group $G$  is  periodic and it follows, like in the finite case, that each element
of $G^{*}$ is of prime order. Moreover, if $p\in \pi(G)$, then each finite $p$-subgroup of $G$  is of order $p$. Therefore
an arbitrary $D(0)$-group  is a $BSP$-group (and in particular a $BCP$-group), as required.
\qed
\enddemo

This Observation implies that results concerning $BSP$-groups (or $BCP$-groups) can be applies to arbitrary $D(0)$-groups.
Therefore  Theorems B6, B7, B3  and B8, dealing with $BSP$-groups (or $BCP$-groups), imply that the corresponding
Theorems 3,4,5 and 6, dealing with $D(0)$-groups, hold.
\qed
\enddemo

We move now to our final section, dealing with arbitrary $D(1)$-groups.

\heading IX. Arbitrary $D(1)$-groups and proofs of Theorems 7,8 and 9\\
\endheading
In this section $G\in D(1)$ still means that $G$ is an arbitrary group belonging to $D(1)$.

And again, if $G$ is a periodic group, then  $\pi(G)$ will denote the set  of all primes dividing the order of some element of $G$.

Recall that  groups $G$ with each  element of $G^{*}$ being of prime power order are called  $CP$-groups. Moreover, a $CP$-group $G$ is called a  $BCP$-group if
for each $p\in \pi(G)$ there exists a positive integer
$u_p$  such that each $p$-element of $G$ is of order $p^i\leq p^{u_p}$, and it is called
a $BSP$-group if  for each $p\in \pi(G)$ there exists a positive integer
$v_p$  such that each finite $p$-subgroup of $G$ is of order $p^j\leq p^{v_p}$.

First we prove the following generalization of Proposition A2.

\proclaim {Proposition C1} If $G\in D(1)$, then the following statements hold.
\roster
\item If $G$ is periodic, then all elements of $G^{*}$ are either of prime order or of prime squared order. In particular,
$G$ is a $BCP$-group.
\item If $G$ is locally finite, then there exists a prime $p \in \pi(G)$ such that a Sylow $p$-subgroup $P$ of $G$
is of cardinality $\geq p^2$, and for all other primes $q\in \pi(G)$, a Sylow $q$-subgroup of $G$ is of the prime order $q$.
In particular, either  $exp(P)=p$ or $exp(P)=p^2$. Furthermore, all elements of $G$ of order $p$ are
deficient and hence they are conjugate in $G$. In the following items, $p$ and $P$ are as defined above.
\item If $G$ is locally finite and $p>2$, then $G$ is a locally solvable group.
\item If $G$ is locally finite and it is not a $p$-group, then $Z(G)=\{1\}$.
\item If $G$ is a locally finite $p$-group, then $G$ is finite.
\item If $G$ is locally finite, then $|\pi(G)|$ is finite and bounded.
\item If $G$ is locally finite, then  every finite subgroup of $G$ of order prime to $p$
has bounded square-free order.
\item If $G$ is locally finite and $P$ is finite, then $G$ is finite.
\item If $G$ is an infinite locally finite group, then  $P$ is infinite.
\endroster
\endproclaim

\demo {Proof}
\roster
\item We can argue as in the proof of (A2(1)).
\item Suppose that $G$ is locally finite and  $p,q$  are distinct primes in $\pi(G)$. Moreover, suppose that    the group $P\in Syl_{p}(G)$
is of cardinality $\geq p^2$
and the group $Q\in Syl_{q}(G)$ is of cardinality $\geq q^2$. Since $G$ is locally finite, there exist finite groups $X$ and $Y$  satisfying $X\leq P$, $Y\leq Q$
and $|X|\geq p^2$, $|Y|\geq q^2$.
If $x\in X$ is a central element of $X$ of order $p$, then $x$ is a deficient element in $G$. Similarly, $Y$ contains a deficient element
of order $q$, contradicting $G\in D(1)$. So there exists  at most one  prime  $p\in \pi(G)$
such that a Sylow $p$-subgroup $P$ of $G$ is of cardinality $\geq p^2$.
If for each $r\in \pi(G)$ a Sylow  $r$-subgroup of $G$ is of prime order $r$, then it follows  by (1) that each element of $G^{*}$ is of prime order
and $G\in D(0)$ by the arguments of Proposition A1, in contradiction
to $G\in D(1)$. Therefore there exists a prime $p\in \pi(G)$ such that a Sylow $p$-subgroup $P$ of $G$ is of cardinality $\geq p^2$
and for each other
prime $q \in \pi(G)$, a Sylow $q$-subgroup of $G$ is of prime order $q$, as required. In particular,
(1) implies that either  $exp(P)=p$ or $exp(P)=p^2$, as required. Finally, suppose that $g\in G$ is of order $p$. Since
$G$ is locally finite
and $P$ is of cardinality $\geq p^2$, there exists a finite subgroup $X$  of $P$ of order $\geq p^2$ and a Sylow $p$-subgroup of $\langle X,g\rangle$ is of finite order
$\geq p^2$. Therefore $C_G(g)>\langle g\rangle$ and $g$ is a deficient element of $G$. Hence  all elements of $G$ of order $p$ are deficient and
since $G\in D(1)$, it follows
that all elements of $G$ of order $p$ are conjugate in $G$, as required.
\item If $G$ is locally finite and  $p>2$, then the order of a Sylow $2$-subgroup of a finite subgroup of $G$
is either $1$ or $2$. Hence $G$ is locally solvable by the Feit-Thompson theorem.
\item If $G$ is locally finite and it is not a $p$-group, and if $x\in Z(G)^{*}$, then there exists  $y\in Z(G)^{*}$ of prime order and hence $G$ contains
an element of order divisible by two primes, in contradiction to (1). Hence $Z(G)=\{1\}$.
\item We may assume that $|G|> p$. By  (2), either $exp(G)=p$ or $exp(G)=p^2$ and all elements of $G$ of order $p$ are conjugate in $G$.
If $p>2$ and $x\in G$ has order $p$, then $x,x^{-1}$ are conjugate . Hence $x^{-1}=x^g$ for some $g\in G$ and $g^2\in C_G(x)$. Since $p>2$, this implies
that $g\in C_G(x)$ and $x^{-1}=x$, a contradiction.

Now suppose that $p=2$. If $exp(G)=2$, then $G$ is abelian and any two non-trivial different elements of $G$ are conjugate, in  contradiction to $|G|>2$.
If $exp(G)=4$ and  $x\in G$ is of order $4$, then $x^2$ is deficient  and therefore $x$ is non-deficient. For any $g\in G$ the subgroup
$\langle x,g\rangle$ is a finite $2$-group, so $1<Z(\langle x,g\rangle)\leq \langle x\rangle$
since $x$ is non-deficient.  Hence $x^2\in Z(\langle x,g\rangle)$, implying that $x^2\in Z(G)$. Since
all elements of $G$ of order $2$ are conjugate in $G$, it follows that $G$ has only one element of order $2$. Since a finite $2$-group with only one
element of order $2$ is either cyclic or generalized quaternion, and since $G$ is a locally finite group with $exp(G)=4$,
it follows that $G$ is either a cyclic group of order $4$ or
a quaternion group of order $8$. In particular, $G$ is finite, as required.
\item  Since by (1) $G$ is a locally finite $CP$-group, it follows by the Main Theorem in [2] that $|\pi(G)|$ is finite
and bounded, as required.
\item Since by (6) $|\pi(G)|$ is finite and bounded, and since by (2) a Sylow $q$-subgroup of $G$ is of prime order for all $q\in \pi(G)\setminus p$,
it follows that a finite subgroup of $G$ of order prime to $p$ has bounded square-free order, as required.
\item Suppose that $P$ is finite. Since by (6) $|\pi(G)|$ is finite and bounded, and since by (2) all Sylow subgroups of $G$ are of finite order, it follows that every finite
subgroup of $G$ is of bounded order. Since $G$ is locally finite, this implies that $G$ is finite, as required.
\item Follows from (8).
\qed
\endroster
\enddemo

\remark {Remark} From now on, if $G$ is a locally finite $D(1)$ group, then $p$ and $P$ will be as defined in (C1(2)).
\endremark

We also need the following proposition dealing with infinite locally finite $D(1)$-groups, which are not $p$-groups.
\proclaim {Proposition C2} Let $G\in D(1)$ be an infinite locally finite group, which is not a  $p$-group.
Then the following statements hold.
\roster
\item Every non-trivial normal subgroup of $G$ is infinite.
\item Every element of $G^{*}$ is of prime order.
\endroster
\endproclaim
\demo {Proof}
\roster
\item Suppose, to the contrary, that $N$ is a finite non-trivial normal subgroup of $G$. Then $M=C_G(N)$ is also normal in $G$ and $G/M$
is a finite group. Hence $M$ is infinite.

Suppose that there exist $a,b\in N$ of prime order satisfying $\langle a\rangle\neq \langle b\rangle$. For every $m\in M\setminus N$ we have
$m\in C_G(a)$ and $m\in C_G(b)$, so $a$ and $b$ are deficient elements of $G$. Since $m$ is not conjugate to $a$ or to $b$, it follows that $m$
is a non-deficient element of $G$. Hence $a\in C_G(m)=\langle m\rangle$ and $b\in C_G(m)=\langle m\rangle$. Since by (C1(1)) $m$ is of prime power
order, it follows that $\langle a\rangle=\langle b\rangle$, a contradiction.

Therefore $N$ is of prime power order, say $q^d$, and it contains only one subgroup of order $q$. Hence $N$ is either cyclic or generalized
quaternion. In any case, there exists a normal subgroup $S$ of $G$ of order $q$. Write $S=\langle s\rangle$. If $q=2$, then $s\in Z(G)$, in contradiction
to (C1(4)). So suppose that $q>2$. Then $C_G(S)$ is a normal subgroup of $G$ of index dividing $q-1$. Hence $C_G(S)$ is an infinite group and by (C1(1)), $C_G(S)$ is a $q$-subgroup of $G$.
If $C_G(S)$ is not locally cyclic, then there exists an element $u\in C_G(S)\setminus S$ of order $q$ and $s,u$ are non-conjugate deficient elements of $G$,
again a contradiction. Therefore $C_G(S)$ is a locally cyclic $q$-group of exponent $\leq q^2$, so it is finite, a
final contradiction.
The proof of (1) is now complete.

\item Suppose, to the contrary, that $x\in G$ satisfies $|x|=p^2$ for some prime $p$. Then by (C1(2)) and (C1(9)) there exists an infinite Sylow
$p$-subgroup $P$ of $G$ containing $x$ and $x^p$ is a deficient element of $G$. Hence $x$ is a non-deficient element of $G$ and for every finite
$p$-subgroup $V$ of $P$ with $x\in V$ and $|V|\geq p^3$, we have $1<Z(V)<\langle x\rangle$ and hence $Z(V)=\langle x^p\rangle$. Since such $V$
could contain any element $y\in P$, it follows that $Z(P)=\langle x^p\rangle$. Since that holds for any element of $P$ of order $p^2$, it follows that
$P/\langle x^p\rangle$ is of exponent $p$.

If $p=2$, then $P/\langle x^2\rangle$ is abelian and $P'\leq \langle x^2\rangle$. In this case $\langle x\rangle$ is normal in $P$, implying that
$|P/C_P(\langle x\rangle)|\leq 2$.
But $x$ is non-deficient in $G$, so $C_P(\langle x\rangle)\leq \langle x\rangle$, implying that $P$ is finite, a contradiction.

So assume that $p>2$.  If  $x\in G$ has order $p$, then by (C1(2)) $x,x^{-1}$ are conjugate .
Hence $x^g=x^{-1}\neq x$ for some $g\in G$ and $g^2\in C_G(x)$. Since $g\notin C_G(x)$, it follows that  $2\mid |g|$ and by (C1(2)) $|g|=2$.
Thus $2\in \pi(G)$. Since by (C1(1)) $G$ is a locally finite $BCP$-group and by (C1(3)) it is locally solvable, it follows by the GH-theorem that $\pi(G)=\{2,p\}$. Since
$p>2$, a Sylow $2$-subgroup of $G$ is of order $2$.
Let $g\in G$  and let $H$ be a finitely generated subgroup of $G$ such that
$x,g\in H$ and $|H|=2p^s$, with $s\geq 3$. Then $H_p$, a Sylow $p$-subgroup of $H$, is normal in $H$ and $|H_p|\geq p^3$.
Clearly $x\in H_p$ and as shown above, $\langle x^p\rangle=Z(H_p)$ and hence
$\langle x^p\rangle$ is normal in $H$. Since $g$ is an arbitrary element of $G$, it follows that  $\langle x^p\rangle$
is normal in $G$, in contradiction to (1).

The proof of the proposition is now complete.
\qed
\endroster
\enddemo

Now we are ready for the proof of the main result of this section.

\proclaim{Theorem 7} If $G\in D(1)$ is a locally finite group, then $G$ is finite.
\endproclaim
\demo{Proof} If $G$ is a $p$-group, then by (C1(5)) $G$ is finite. So suppose that $G$ is not a $p$-group and suppose
that $G$ is infinite. Our aim is to reach a contradiction.

By (C2(2)) every element of $G^{*}$ is of prime order. By (C1(2)),  $P$ is
of cardinality $\geq p^2$ and of exponent $p$, while for each prime $q\in \pi(G)\setminus p$ a Sylow $q$-subgroup of $G$ is of prime order $q$.
Moreover, the following statements hold: by (C1(6)) $|\pi(G)|$ is finite and bounded, by (C1(9)) $P$ is infinite and by (C1(2))  all elements of $G$ of order $p$ are conjugate.

By Theorem 5 of [2], if $H$ is a locally finite simple group with all non-trivial elements of prime order, then $H$ is finite. Therefore
$G$ is non-simple, and there exists a proper  non-trivial normal subgroup $N$ of $G$. If $p\notin \pi(N)$, then by (C1(7)) every finite subgroup of $N$
has a finite  bounded order and since $N$ is locally finite, it is finite, in contradiction to (C2(1)). Hence $p\in \pi(N)$ and by (C1(2)), $N$ contains all elements of $G$ of
order $p$. Hence $P\leq N$ and $p\notin \pi(G/N)$. Moreover, since $G/N$ is locally finite with $|\pi(G/N)|$ finite and bounded and with all Sylow subgroups of prime order, it follows that a finite subgroup of $G/N$ has a bounded order and hence $G/N$ is finite. For any $a\in N^{*}$ we have $C_G(a)\leq N$, since otherwise if $b\in C_G(a)\setminus N$,
then $a,b$ are non-conjugate commuting elements of $G$ of prime order and  only one of them may be deficient, so $\langle a\rangle=\langle b\rangle$ and $b\in N$,
a contradiction. Therefore $G$ is a locally finite Frobenius group with the kernel $N$, implying that $N$ is nilpotent (see Theorem 1.J.2 in [9]).
Hence $Z(N)\neq \{1\}$ and if $x\neq 1$ belongs to $Z(N)$, then $x\in C_G(P)$. This implies that $|x|=p$ and $N\leq  C_G(x)$.
Since $G/N$ is finite, it follows that the conjugacy class $x^G$ of $x$ in $G$
is finite. But $P^{*}$ is contained in $x^G$, so $P$ is finite, a final contradiction.
The proof of Theorem 7 is now complete.
\qed
\enddemo

Theorem 7 yields the following two corollaries.
\proclaim{Theorem 8} If $G$ is a periodic locally graded $D(1)$-group, then $G$ is finite.
\endproclaim
\demo{Proof} By (C1(1)), a periodic $D(1)$-group is a $BCP$-group. Hence $G$ is a locally graded $BCP$-group, so  by Theorem B4 $G$ is a
locally finite $D(1)$-group.
It follows then by Theorem 7 that $G$ is finite, as required.
\qed
\enddemo
\proclaim{Theorem 9} If $G$ is a residually finite $D(1)$-group, then $G$ is finite.
\endproclaim
\demo{Proof} First we show that $G$ is a periodic group. Suppose, to the contrary, that $x\in G$ is of infinite order.
Since $G$ is residually finite, there exists a normal subgroup $N$ of $G$ of finite index $n>1$ such that $x \notin N$ and $x^n\in N$.
Moreover, there exists a normal subgroup $M$ of $G$ of finite index $m>1$, such that $x^n\notin M$. Then $x^m\in M$ and
$x^m$ is not conjugate in $G$ to $x^n$.
Since $\langle x^n\rangle<\langle x\rangle\leq C_G(x^n)$ and $\langle x^m\rangle<\langle x\rangle\leq C_G(x^m)$, it follows that
$x^n$ and $x^m$ are non-conjugate deficient elements of $G$, in contradiction to $G\in D(1)$.
Therefore $G$ is a periodic $D(1)$-group.

By (C1(1)), a periodic $D(1)$-group is a $BCP$-group. Hence $G$ is a residually finite $BCP$-group, so  by Theorem B3 $G$ is a
locally finite group.
It follows then by Theorem 7 that $G$ is finite, as required.
\qed
\enddemo

\enddocument